\documentclass[reqno]{amsart}

\usepackage{amsfonts}
\usepackage{amsmath}
\usepackage{amscd}
\usepackage{amssymb}
\usepackage{diagrams}
\usepackage{verbatim}

\newcommand{\Z}{\mathbb{Z}}

\newcommand{\N}{\mathbb{N}}
\newcommand{\F}[1][p]{\mathbb{F}_{#1}}

\newcommand{\CP}[1][\infty]{\mathbb{C}P^{#1}}   
\newcommand{\E}[1][E]{\mathbf{#1}}              
\newcommand{\EHAT}[1][n]{\E[\widehat{E(#1)}]}   
\newcommand{\Ehat}[1][n]{\widehat{E(#1)}}       
\def\K{\mathbf{K(n)}}                           

\newcommand{\Ch}{\mathit{Ch}}                   
\newcommand{\ECh}[1][E]{\Ch_{#1}}               
\newcommand{\EhatCh}[1][n]{\Ch_{\widehat{E(#1)}}}

\newcommand{\A}{\mathcal{A}}                    
\newcommand{\An}[1][n]{\A^{(#1)}}               
\newcommand{\Ae}[1][E]{\A_{#1}}                 
\newcommand{\C}{\mathcal{C}}                    
\newcommand{\CR}[1][R]{\mathcal{C}(#1)}         
\newcommand{\CE}[1][E]{\mathcal{C}_{#1}}        

\newcommand{\Res}{\operatorname{Res}}               
\newcommand{\Tr}{\operatorname{Tr}}                 
\newcommand{\Hom}{\operatorname{Hom}}               
\newcommand{\var}{\operatorname{var}}               
\newcommand{\colim}{\operatornamewithlimits{colim}} 
\newcommand{\Vark}[1]{\var(#1)}			    
\newcommand{\tensor}{\displaystyle\mathop{\otimes}} 

\newcommand{\lra}{\longrightarrow}
\newcommand{\w}{\mathbf{w}}
\newcommand{\barb}{\bar{b}}
\newcommand{\onto}{\twoheadrightarrow}
\newcommand{\Fiso}{F}

\newtheorem{thm}{Theorem}[section]

\newtheorem{lemma}[thm]{Lemma}
\newtheorem{prop}[thm]{Proposition}
\newtheorem{example}[thm]{Example}
\newtheorem{defn}[thm]{Definition}

\newtheorem{remark}[thm]{Remark}
\newtheorem{property}[thm]{Property}

\begin{document}

\title[Chromatic characteristic classes]
{Chromatic characteristic classes\\ in ordinary group cohomology}

\author[D.~J.~Green]{David J.\ Green}
\author[J.~R.~Hunton]{John R.\ Hunton}
\author[B.~Schuster]{Bj\"orn Schuster}

\date{20 Aug 2001}

\address{Department of Mathematics, University of Wuppertal, Gau{\ss}str.~20,
D-42097 Wuppertal, Germany}
\email{green@math.uni-wuppertal.de, schuster@math.uni-wuppertal.de}

\address{Department of Mathematics and Computer Science, University Road,
Leicester, LE1 7RH, England}
\email{J.Hunton@mcs.le.ac.uk}

\begin{abstract}
We study a family of subrings, indexed by the natural numbers, of the mod $p$
cohomology of a finite group $G$. These subrings are based on a family of
$v_n$-periodic complex oriented cohomology theories and are constructed as
rings of generalised characteristic classes. We identify the varieties
associated to these subrings in terms of colimits over categories of
elementary abelian subgroups of $G$, naturally interpolating between the work 
of Quillen on $\Vark{H^*(BG)}$, the variety of the whole cohomology ring, and 
that of Green and Leary on the variety of the Chern subring, $\Vark{Ch(G)}$. 
Our subrings give rise to a \lq chromatic\rq\ (co)filtration, which has both
topological and algebraic definitions, of $\Vark{H^*(BG)}$
whose final quotient is the variety $\Vark{Ch(G)}$.
\end{abstract}

\maketitle

\section*{Introduction}
\noindent
This paper develops a structure within the mod $p$ cohomology ring
$H^*(BG)=H^*(BG;\F[p])$ of a finite group $G$ which we show to arise
both topologically and algebraically. In general, the determination
and description of the ring $H^*(BG)$ is a hard problem and we build
upon two separate approaches, due respectively to  Quillen \cite{Quillen}
and to Thomas \cite{Thomas}, incorporating also ideas from the chromatic
point of view in homotopy theory (in particular \cite{HKR}) and the theory
of coalgebraic rings (Hopf rings) and the homology of infinite loop spaces
\cite{BH,hun,HT1,RW,wkn}. We utilise the work of Leary and the first author
\cite{VCh} on categories associated to subrings of $H^*(BG)$, which first
hinted at the possibility of a theorem such as our result \eqref{theorem:main}
below.

\medskip

Recall that Quillen showed \cite{Quillen} that the ring $H^*(BG)$ could
be described up to \Fiso-isomorphism in terms of its elementary
abelian subgroup structure. For its most elegant formulation, suppose
$k$ is an algebraically closed field of characteristic $p$ and for an
$\F[p]$ algebra $R$ denote by $\Vark{R}$ the variety of algebra morphisms
from $R$ to $k$, topologised with the Zariski topology. Quillen's result
describes $\Vark{H^*(BG)}$ in terms of a colimit $\colim_{\A}\Vark{H^*(BV)}$
over a certain category $\A$ of elementary abelian subgroups $V$ of $G$.
In particular recall that, for such $V$, the variety $\Vark{H^*(BV)}$ is
isomorphic to $V\tensor_{\F[p]}k$ where $V$ is viewed as an $\F[p]$
vector space.

On the other hand, Thomas has shown \cite{Thomas} that a frequently
significant as well as computable subring of $H^*(BG)$ is the
{\em Chern subring}, $Ch(G)$, essentially that part of the cohomology
ring given most immediately by the complex representation theory of $G$
(specifically, it is the subring generated by all Chern classes of irreducible
representations). Leary and the first author \cite{VCh} have completed this
picture by describing $\Vark{Ch(G)}$ as $\colim_{\A^{(1)}}\Vark{H^*(BV)}$,
the colimit over a further category $\A^{(1)}$ of elementary abelians.
The categories $\A$ and $\A^{(1)}$ have the same objects but in general
they have different morphisms and the natural surjection
$\Vark{H^*(BG)}\onto \Vark{Ch(G)}$ is not usually a homeomorphism.

\medskip
However, Atiyah's theorem \cite{Atiyah} linking the (completed) representation
theory of $G$ with $K^*(BG)$, the complex $K$-theory of $BG$, allows a
reinterpretation of $Ch(G)$ in terms of (unstable) maps from $BG$ to the
spaces in the $\Omega$-spectrum for $K$-theory: the homology images of such
maps can be identified with Chern classes, the generators of the Chern subring.
See section~\ref{section:ecc} for details. Replacing $K$-theory by any other
representable generalised cohomology theory $E^*(-)$, this construction
allows for the definition of `$E$-type characteristic classes' in $H^*(BG)$
based on the spaces in the $\Omega$-spectrum for $E$. These generate a
subring, $Ch_E(G)$ say, of $H^*(BG)$. In many instances (including, by the
work of Bendersky and the second author \cite{BH}, those when $E$ is Landweber
exact \cite{Landweber}) there are inclusions of rings
$$Ch(G)\subset Ch_E(G)\subset H^*(BG).$$

Our main results concern a family of such subrings, defined by a family
of cohomology theories $E$. Specifically, we are interested in the subrings
given by a family of $v_n$-periodic cohomology theories such as the
Johnson-Wilson \cite{JW} theories $E(n)$, $n\in\N$. However, to obtain our best
results we concentrate on the subrings given by their $I_n$-adic completions
$\widehat{E(n)}$, as given by Baker and W\"urgler \cite{bwa}, and these
subrings we denote $\EhatCh(G)$. The use of complete theories has the
advantage of enabling us to apply the results of Hopkins, Kuhn and Ravenel
\cite{HKR}. In section~\ref{section:thecat} we define certain categories
$\An$ of elementary abelian subgroups of $G$ allowing us to prove

\begin{thm}
\label{theorem:main}
Let $G$ be a finite group and $k$ an algebraically closed field of
characteristic $p$. Then there is a homeomorphism of varieties
$$\colim_{\An} \Vark{H^*(BV)} \longrightarrow \Vark{\EhatCh(G)}\,.$$
\end{thm}

\smallskip
The categories $\An$ appeared in \cite{VCh} where, as mentioned, the first
case $\A^{(1)}$ was proved to represent the classical Chern subring. It was
also suggested that there may be a possible link between the $\An$ and certain
$v_n$-periodic spectra; this theorem provides such a link and supplies a
topological construction for subrings of $H^*(BG)$ associated to the $\An$.

\medskip
\noindent
The same result almost certainly holds on replacing the spectra
$\widehat{E(n)}$ by Morava $E$-theory or in fact by any Landweber exact
spectrum satisfying the main hypotheses of \cite{HKR}; the recent work
of Baker and Lazarev \cite{BL} supplies the technical constructions needed
for adapting the proof given below for $\widehat{E(n)}$ (the essential element
being the analogue of the Baker-W\"urgler tower \cite{bwb}). It is less clear
that the analagous result for $E(n)$ will hold as these spectra fail the
completion hypotheses of \cite{HKR}.

\medskip
A number of basic results on the $\EhatCh(G)$ can be proved.

\begin{thm}\label{theorem:basics}
\begin{enumerate}
\item The varieties $\Vark{\EhatCh(G)}$ form a (co)filtered space
$$\Vark{H^*(BG)}\onto\cdots\onto \Vark{\Ch_{\widehat{E(n+1)}}(G)}\onto
\Vark{\Ch_{\widehat{E(n)}}(G)}\onto\cdots\onto k$$

\item For each $n\geqslant 1$ there is a finite group $G$ for which
$\Vark{\EhatCh(G)}$ is distinct from $\Vark{\Ch_{\widehat{E(n+1)}}(G)}$,
and hence $\EhatCh(G)$ from $\Ch_{\widehat{E(n+1)}}(G)$;

\item For $n$ at least the $p$-rank of $G$, $\Vark{\EhatCh(G)}=\Vark{H^*(BG)}$;

\item $\Vark{\Ch_{\widehat{E(1)}}(G)}=\Vark{Ch(G)}$, the Chern subring of
$H^*(BG)$;

\item For each $n$ there is an inclusion of rings $Ch(G)\subset\EhatCh(G)$.
\end{enumerate}
\end{thm}

Given Theorem~\ref{theorem:main}, the statements (1) and (3) follow
immediately from the definition of the $\An$ while (2) and (4) follow
from results of \cite{VCh}. A more general version of statement (5) is
proved in section~\ref{section:VCh}.

For simplicity we shall assume, at least after section~\ref{section:ecc},
that all our rings $R$, such as $Ch_E(G)$, are defined as subrings of
the even degree part of $H^*(BG)$, and hence are commutative. In fact,
as noted in \cite{VCh}, this is not restrictive: if $p\not=2$ the odd
degree elements of $H^*(BG)$ are nilpotent and any homomorphism
$H^{\rm even}(BG)\supset R\buildrel f\over\rightarrow k$ extends
trivially for any odd degree element; if $p=2$ such a map would extend
uniquely on an odd degree element $x$ to the unique square root of $f(x^2)$.

\medskip
The rest of the paper is organised as follows. We begin by defining our
$E$-type characteristic classes in section~\ref{section:ecc}, linking this
construction to that of the classical Chern subring in the case $E=K$. We
also introduce here some of the basic notation and language of coalgebraic
algebra \cite{HT1,RW} in order to discuss the homology of spaces in the
$\Omega$-spectrum for $E$.

In section~\ref{section:VCh} we recall and extend the main result of \cite{VCh}
and prove that, at least for $E$ Landweber exact, $\Vark{Ch_E(G)}$ has a
description as a colimit of $\Vark{H^*(BV)}$ as $V$ runs over some category,
$\A_E$ say, of elementary abelian subgroups $V$ of $G$. Here we also
demonstrate that for such $E$ there is an inclusion $Ch(G)\subset Ch_E(G)$.

Most of the rest of the paper is devoted to the proof of the
Theorem~\ref{theorem:main}, {\em i.e}, to the identification of the
category $\A_E$. In section~\ref{section:thecat} we introduce the cateogry
$\An$ together with an intermediate category $\CE$. We show that if $E$ both
is Landweber exact and satisfies the assumptions of the character isomorphism
of \cite{HKR} then it follows fairly easily that $\CE=\An$. The harder part
is to show that $\A_E=\CE$ and we prove that this will follow if certain
properties of the `unstable Hurewicz homomorphism'
$$\mathcal{H}\colon E^*(X)=[X,\E_*]\rightarrow\Hom(H^*(\E_*),H^*(X))$$
can be established. Our most general result, which for $E=\Ehat$ specialises to
Theorem~\ref{theorem:main}, is stated as Theorem~\ref{theorem:maingen} and the
complete set of hypotheses needed on a spectrum $E$ for its application are
set out in \eqref{cond}.

Proof that these hypotheses are satisfied by $\Ehat$ is given in
section~\ref{section:theproof} which examines in detail the homology
of the spaces in its $\Omega$-spectrum, building on Baker and W\"urgler's
description \cite{bwb} of $\widehat{E(n)}$ as a homotopy limit and Wilson's
calculation of the Hopf ring for Morava $K$-theory \cite{wkn}.

We conclude with some example computations in section~\ref{section:examples}.

\medskip
The construction of the subrings $Ch_E(G)$ is firmly based in the world of
unstable homotopy through the use of the infinite loop spaces representing
$E$-cohomology. We rely extensively on the techniques of coalgebraic algebra
for our proofs and computations; see \cite{HT1,RW} for introductory material
on this field. We also rely on the work of Wilson \cite{wkn} as the major
computational input, but we note one more recent aspect of our work. The
results of \cite{HKR} needed in section~\ref{section:thecat} demand we work
with completed spectra $E$ and consequently with the spaces in their
$\Omega$-spectra. These are truly huge spaces whose homologies have only
recently become accessible; see \cite{BH,hun,HT2} for the methods of
handling such spaces and for the results we use below.

\medskip
\smallskip\noindent{\bf Acknowledgements} We are grateful for the
suggestions and advice of Andy Baker, Dave Benson, John Greenlees, Ian
Leary, Erich Ossa and
Steve Wilson. The second author thanks the Bergische Universit\"at Wuppertal
for financial support and the University of Leicester for sabbatical leave
while much of this work was completed.

\medskip

\section{$E$-type characteristic classes}
\label{section:ecc}
In this section we define the class of subrings considered and show that
the classical Chern subring arises as a special case. Here and in the next
section we can take $G$ to be any compact Lie group; later we shall need to
restrict to finite groups. We begin with a very general construction.

\begin{defn}\label{general}
Suppose $\{X_i\}$ is a family of spaces and $\mathcal{F}G$ is a set of
maps of the form $f\colon BG\rightarrow X_i$. Let $Ch_{\mathcal{F}G}$
be the subring of $H^*(BG)$ generated by all elements of the form $f^*(x)$
as $f$ runs over the elements of $\mathcal{F}G$ and $x$ over the homogeneous
elements of all the $H^*(X_i)$.
\end{defn}

\begin{example}\label{standardchern}{\em
Take $\{X_i\}$ to be the set containing just one space, $BU$. Let
$\mathcal{F}$ be the set of maps $BG\rightarrow BU$ given by irreducible
representations $G\rightarrow U(n)$ on embedding $U(n)$ in the infinite
unitary group $U$ and taking classifying space. The subring in this instance
is the classical Chern subring $Ch(G)$ \cite{Thomas}.}
\end{example}

\begin{remark}\label{RGchern}{\em
It is not hard to show that the Chern subring, as given in the last
example, can equivalently be defined by taking $\mathcal{F}G$ as the
maps given by all representations of $G$, or even as the maps given
by all elements of the representation ring $R(G)$. Clearly enlarging
$\mathcal{F}G$ will only if anything enlarge the subring of $H^*(BG)$
it defines, so it suffices to show that elements of $R(G)$ give no new
elements than are already in $Ch(G)$ as defined in~\eqref{standardchern}.
For example, if $\rho_1\colon G\rightarrow U(n)$ and
$\rho_2\colon G\rightarrow U(m)$ are two representations, the H-space
structure on $BU$ induced by $U(n)\times U(m)\rightarrow U(n+m)$ gives,
in cohomology, the commutative diagram
\medskip
\begin{equation*}
\label{eqn:RGchern}
\begin{CD}
H^*(BU)\otimes H^*(BU) @>{B\rho_1^*\otimes B\rho_2^*}>> H^*(BG)\otimes
H^*(BG) \\
@AAA @VVV \\
H^*(BU) @>{B(\rho_1+\rho_2)^*}>> H^*(BG)
\end{CD}
\end{equation*}

\medskip\noindent
sending an element $x\in H^*(BU)$ to $\sum B\rho_1^*(x')\cdot B\rho_2^*(x'')$
in $H^*(BG)$ where $\sum x'\otimes x''$ denotes the image of $x$ in
$H^*(BU)\otimes H^*(BU)$. Thus the representation $\rho_1+\rho_2$ fails to
give anything in $H^*(BG)$ that was not already in the subring defined by
$\rho_1$ and $\rho_2$. The rest of the proof is similar.}
\end{remark}

To set up notation that will be useful later, recall that $H^*(BU)$ is
the polynomial ring on generators $c_1,c_2,\ldots$, the universal Chern
classes, where $c_i\in H^{2i}(BU)$. The pull-back $B\rho^*(c_i)$ of a
representation $\rho\colon G\rightarrow U$ is the called the $i^{\rm th}$
Chern class of $\rho$, also written $c_i(\rho)$, and is an element of
$H^{2i}(BG)$.

\medskip
For the remainder of this paper we assume $E^*(-)$ is a representable
cohomology theory on the category of CW complexes with coefficients $E_*$
concentrated in even dimensions. In particular, there is a representing
$\Omega$-spectrum with spaces $\E_r$, $r\in\Z$, and, for all CW complexes $X$,
natural equivalences between $E^r(X)$ and the set of homotopy classes of maps
from $X$ to $\E_r$. The following is the main definition of this paper.

\begin{defn}\label{maindef}
Let $E$ be as just described and let $G$ be a compact Lie group.
Let $\mathcal{F}G$ be the set of all homotopy classes of maps
$BG\rightarrow\E_{2r}$, allowing all $r\in\Z$. The {\em $E$-Chern subring},
$Ch_E(G)$, is defined as the subring $Ch_{\mathcal{F}G}$ of $H^*(BG)$.
Also, define an {\em $E$-type characteristic class\/} of~$G$ to be an
element $f^*(x) \in H^*(BG)$ for some $f\colon BG\rightarrow\E_{2r}$ and
some homogeneous $x\in H^*(\E_{2r})$.
\end{defn}

\begin{remark}\label{conncmpt}{\em
Note that it does not matter in this definition whether we take the whole
spaces $\E_{2r}$ or just the connected components of their basepoints, usually
denoted $\E_{2r}'$. This follows since all components of the $\E_{2r}$ are the
same, with $\E_{2r}=\pi_0(\E_{2r})\times\E_{2r}'=E^{2r}\times\E_{2r}'$.}
\end{remark}

The following result, building on Remark~\ref{RGchern}, shows that this
definition includes that of the classical Chern subring.

\begin{prop}\label{prop:OKreally}
Let~$G$ be a compact Lie group. Then the $K$-Chern subring $\ECh[K](G)$ is
equal to $\Ch(G)$, the classical Chern subring.
\end{prop}

\begin{proof}
As any element of the representation ring $R(G)$ gives rise to a homotopy
class of maps $BG\rightarrow BU$, it is immediate by \eqref{RGchern} and
\eqref{conncmpt} that every element of $Ch(G)$ is also an element of
$Ch_K(G)$ since the only even graded spaces in the $\Omega$-spectrum
for $K$-theory are of the homotopy type of $\Z\times BU$. However, not
every class of map $BG\rightarrow BU$  necessarily arises as the image
under $B(-)$ of something from the representation ring. Thus we must check
further that $Ch_K(G)\subset Ch(G)$.

Recall from the work of Atiyah and Segal~\cite{AS} (or Atiyah~\cite{Atiyah}
for the case of finite groups) that the homomorphism
$\alpha\colon R(G)\rightarrow K^0(BG)$ given by sending
the virtual representation $\rho$ to the corresponding map
$BG\rightarrow\Z\times BU$ is continuous with respect to the
$I$-adic topology on $R(G)$, where $I$ is the augmentation ideal,
and the filtration topology on $K^0(BG)$ given by
$$F_mK^0(BG)=\ker(K^0(BG)\rightarrow K^0(BG^{m-1})) \, ,$$
where $BG^{m-1}$ is the $m-1$ skeleton of $BG$. Moreover, $K^0(BG)$
is a complete ring and $\alpha$ is an isomorphism after $I$-adic
completion of $R(G)$.

Note that the degree $m$ cohomology of $BG$ is determined by its
$(m+1)$-skeleton $BG^{m+1}$, and so in particular the $r^{\rm th}$
Chern class $c_r(g)$ of a map $g\colon BG\rightarrow BU$ is determined
by its restriction to $BG^{2r+1}\rightarrow BU$. It follows that $c_r(g)$
is zero for all $1\leqslant r\leqslant n-1$ if $g\in F_{2n}K^0(BG)$.

Now consider a map $f\colon BG\rightarrow BU$ in $\mathcal{F}G$.
As an element of $K^0(BG)$ it is represented by a Cauchy sequence
$\rho_1,\rho_2,\ldots$ in $R(G)$. By the Whitney sum formula, the
total Chern class of $f$ satisfies $c(f) = c(\rho_i) c(f - B\rho_i)$.
Hence if $f - B\rho_i\in F_{2n}K^0(BG)$, then $c_r(f)=c_r(\rho_i)$
for $1\leqslant r\leqslant n-1$. Thus, for each $r\geqslant 0$, the
sequence $c_r(\rho_1), c_r(\rho_2),\ldots$ in $H^{2r}(BG)$ is eventually
constant at $c_r(f) = f^*(c_r)$.

As any homogeneous $x\in H^*(BU)$ can involve only finitely many of the
$c_r$, each element $f^*(x)\in H^*(BG)$ can be realised in the form
${B\rho}^*(x)$ for some virtual representation $\rho \in R(G)$. This
completes the proof.
\end{proof}

\begin{remark}{\em
These definitions of $E$-type characteristic classes can also be made in
the integral or $p$-local rings $H^*(BG;\Z)$ and $H^*(BG;\Z_{(p)})$. If
$H^*(\E_{2*};\Z)$ is torsion free (as is the case when $E$ is Landweber
exact, \cite{BH}), then the (mod-$p$) $E$-type characteristic classes are
just the mod-$p$ reductions of the integral $E$-type characteristic classes.
In general however this need not be the case. As $H^*(BU;\Z)$ is torsion free,
the integral analogue of \eqref{prop:OKreally} holds.}
\end{remark}

\begin{remark}{\em
The reader may wonder as to the restriction in Definition~\ref{maindef} to
even graded spaces $\E_{2r}$ in the $\Omega$-spectrum for $E$. In the case
of complex $K$-theory and finite groups $G$ this restriction is vacuous:
there are no non-trivial maps from $BG$ to the odd graded spaces \cite{Atiyah}.
For more general compact groups Proposition~\ref{prop:OKreally} would however
fail without this hypothesis and so \eqref{maindef} would not have been the
correct extension of the concept of Chern subring.
Below we shall restrict to finite $G$ and Landweber exact $E$;
in light of Kriz's examples such as presented in \cite{Kriz}, it is
conceivable that even here relaxing Definition~\ref{maindef} to include all
spaces could enlarge the subring it defines, but the effect remains unclear.
The theorems below are proved considering the even spaces alone.
}\end{remark}

\begin{remark}{\em
There is of course no {\em a priori\/} reason why the subring $Ch_E(G)$
defined by \eqref{maindef} should lie in $H^{\rm even}(BG)$, since there
may well be odd dimensional elements of $H^*(\E_{2r})$ which pass to
non-trivial elements of $H^*(BG)$. However, we shall need from the main
result of the next section, Theorem~\ref{theorem:large}, and onwards that
the spectrum $E$ considered is Landweber exact. We can thus appeal to the
work of \cite{HT2} from which it can easily be deduced that whenever $E$
is Landweber exact $H^*(\E_{2r})$ has no odd dimensional elements. }
\end{remark}

\section{Representations by categories of elementary abelians}
\label{section:VCh}

Quillen proved in~\cite{Quillen} that the variety of~$H^*(BG)$ can be built
from the elementary abelian $p$-subgroups using conjugacy information. This
result was extended in~\cite{VCh} to cover the varieties of certain large
subrings of~$H^*(BG)$ as well. The object of this section is to show that for
a wide range of theories $E$, the varieties of the $E$-Chern subrings $Ch_E(G)$
have a similar description.  We start by recalling the definitions and the main
result of \cite{VCh}, in a slightly sharpened form.

\begin{defn}
For a compact Lie group~$G$ and homogeneously generated subring $R$ of
$H^*(BG)$, define the category $\CR$ as that with objects the elementary
abelian $p$-subgroups $V$ of $G$, and morphisms from $W$ to $V$ to be the
set of injective group homomorphisms $f \colon W \rightarrow V$ satisfying
\begin{equation}\label{CRdef}
f^* \Res_{V} (x) = \Res_{W} (x)
\end{equation}
for every homogeneous $x \in R$.
\end{defn}

\begin{remark}{\em
This is a variation on the definition given in \cite{VCh} where the
condition on morphisms $f\colon W\rightarrow V$ was weakened to
satisfying equation~\eqref{CRdef} modulo the nilradical.
However, for subrings $R$ that are closed under the action of the Steenrod
algebra, the two definitions are equivalent.
To see this, suppose that $f$ is a morphism in the ``modulo nilradical''
version of $\CR$ with the property that,
for some $x\in R$, the class $w=f^*\Res_{V}(x) - \Res_{W} (x)$ is
a non-zero nilpotent in $H^*(BW)$. For any non-zero homogeneous element
of~$H^*(BW)$ there is an operation $\theta$ in the mod $p$ Steenrod algebra
such that $\theta(w)$ is  non-nilpotent: this is immediate for $p=2$; for $p$
odd, see Lemma~2.6.5 of~\cite{Schwartz}. Naturality of~$\theta$ implies that
$$f^* \Res_{V} (\theta x) - \Res_{W} (\theta x) = \theta(w) \, .$$
As $R$ is closed under the Steenrod algebra, $\theta(x)$ lies in $R$,
contradicting $f$ being a morphism in $\CR$.
}\end{remark}

\begin{defn}\label{def:large}
Let~$G$ be a compact Lie group.  A virtual representation $\rho$~of $G$
is called \emph{$p$-regular} if it has positive virtual dimension and restricts
to every elementary abelian $p$-subgroup as a direct sum of copies of the
regular representation.
A subring $R$~of $H^*(BG)$ is called \emph{large} if it contains the Chern
classes of some $p$-regular representation.
\end{defn}

\begin{thm}
\label{theorem:VCh} \cite[6.1]{VCh}
Let $G$ be a compact Lie group and $R$ a homogenously generated subring
of~$H^*(BG)$ which is both large and closed under the Steenrod algebra.
Then the natural map $R \rightarrow \lim_{V \in \CR} H^*(BV)$ induces a
homeomorphism
\begin{equation}
\colim_{V \in \CR} \Vark{H^*(BV)} \rightarrow \Vark{R} \, . 
\tag*{\qedsymbol}
\end{equation}
\end{thm}

The problem, however, is to identify the category $\CR$ for a given $R$.

\begin{example}{\em
The Quillen category~$\A = \A(G)$ is the category with objects the
elementary abelian $p$-subgroups of~$G$, and morphisms generated by
inclusion and conjugation.  It is shown in \cite[9.2]{VCh} that $\A$~is
the category~$\CR[H^*(BG)]$, thus recovering one of the main results
of~\cite{Quillen}:
\[ \colim_{V \in \A} \Vark{H^*(BV)} \cong \Vark{H^*(BG)} \, . \]
}\end{example}

\begin{example}\label{example:A1}{\em
It is shown in \cite{VCh} that the category $\A^{(1)}$ of elementary
abelian $p$-subgroups of $G$ and morphisms the injective homomorphisms
$f\colon W\rightarrow V$ which take each element $w\in W$ to a conjugate
of itself is the category $\CR[Ch(G)]$ associated to the Chern subring.
}\end{example}

We come to the main result of this section. Recall that the class of Landweber
exact spectra \cite{Landweber} includes examples such as $BP$-theory, complex
cobordism, the Johnson-Wilson theories $E(n)$ and their $I_n$-adic completions
$\widehat{E(n)}$ as well as Morava $E$-theory, complex $K$-theory, various
forms of elliptic cohomology \cite{LRS} and their completions (in particular
the completions $(Ell)^{\widehat{\quad}}_\wp$ of Baker \cite{bak}).

\begin{thm}
\label{theorem:large}
Let $E$ be a Landweber exact spectrum and let $G$ be a compact Lie group.
Then $Ch_E(G)$ is closed under the Steenrod algebra and is large in the
sense  of {\em \eqref{def:large}.} Hence there is a category $\CR[\ECh(G)]$,
which we shall write as $\Ae(G)$ or just $\Ae$, such that
$$\Vark{\ECh(G)} \cong \colim_{V \in \Ae} \Vark{H^*(BV)} \, . $$
\end{thm}

\begin{proof}
Closure of $Ch_E(G)$ under the Steenrod algebra is immediate from its
definition and it
suffices to show that it is also large. This follows from the next result
and the fact
that the Chern subring $Ch(G)$ is large \cite{VCh}.
\end{proof}

\begin{prop}\label{prop:inclusion}
For $E$ Landweber exact the homomorphism
$(c_1^E)^* \colon H^*(\E_2') \rightarrow H^*(BU)$ is surjective for all
$n \geqslant 1$. Hence for such $E$ there is an inclusion of subrings
$Ch(G)\subset Ch_E(G)$ for all compact Lie groups $G$ and in particular
$Ch_E(G)$ contains all Chern classes of complex representations of $G$.
\end{prop}

\noindent Here $c_1^E$ is the first $E$-theory Chern class: recall that
$$E^*(BU)=E^* [\![c^E_1, c^E_2, \ldots ]\!]\ \mbox{ with }\ c^E_i\in
E^{2i}(BU) \, ; $$
thus $c^E_1 \colon BU \rightarrow \E_2'$ is an element of $\widetilde{E}^2(BU)$
where, as before, $\E_2'$ means the base point component of $\E_2$.

We need to introduce further notation. Recall that any Landweber exact
spectrum $E$ is certainly complex oriented. For any complex oriented spectrum,
denote by $x^E\colon\CP\rightarrow \E_2'$ in $\widetilde{E}^2(\CP)$ the
orientation map. The definition of $c^E_1$ is not independent of this and
in fact $c_1^E\colon BU\rightarrow\E_2'$ restricts to the complex orientation
$x^E$ on $\CP$ \cite[Part II, 4.3({\em ii})]{Adams}. We assume that the
orientation $x^K$ has been chosen compatibly so that the diagram
\begin{equation*}
\begin{diagram}[height=2em,width=2em,tight]
   BU &                & \rTo^{c_1^E} &             & \E_2' \\
      & \luTo_{x^K}    &              & \ruTo_{x^E} &       \\
      &                &     \CP      &             &       \\
\end{diagram}
\end{equation*}
commutes.

Take basis elements $\beta_n$, $n\in\N$, of $H_*(\CP)$ with
$\beta_n\in H_{2n}(\CP)$, dual to the basis of $H^*(\CP)$ given by the
powers of the orientation $x^H$. For a complex oriented theory $E$
define the elements $b_n^E\in H_n(\E_2')$ by $b_n^E=(x^E)_*(\beta_n)$.

Finally, as $E^r(X)$ is an abelian group for any space $X$, the space
$\E_r$ has an H-space product $\E_r\times\E_r\rightarrow\E_r$ which leads,
in mod $p$ homology, to a product
$$*\colon H_*(\E_r)\otimes H_*(\E_r)\rightarrow H_*(\E_r).$$

\begin{lemma}
Let $E$ be any complex oriented spectrum and let $K$ denote complex $K$-theory.
Then $(c_1^E)_*\left((b_1^K)^{*r}\right) =  (b_1^E)^{*r}$ for all~$r$.
\end{lemma}

\begin{proof}
The diagram above, in homology, says $(c_1^E)_*(b_1^K)=b_1^E$. The
lemma follows by noting that $c_1^E\in E^2(BU)$ is a primitive element
in the Hopf algebra $E^*(BU)$, and hence represents an (unstable) additive
cohomology operation $K^0(-)\rightarrow E^2(-)$. Thus the map $c_1^E$
commutes with $*$-products and
$$(c_1^E)_*\left((b_1^K)^{*r}\right) =
\left((c_1^E)_*(b_1^K)\right)^{*r} =(b_1^E)^{*r}$$
for all~$r$.
\end{proof}

\begin{proof}[Proof of Proposition~\ref{prop:inclusion}]
By the work of Bendersky and the second author~\cite{BH},
$H_*(\E_{2r})$ is polynomial (under the $*$-product) for any
Landweber exact $E$ and $b_1^E\ne 0$. It follows that
$(c_1^E)_*({b_1^K}^{*r}) = ({b_1^E})^{*r}$ is never zero.
But $c^K_r \in H^{2r}(BU)$ is dual to $(b^K_1)^{*r} \in H_{2r}(BU)$,
with respect to the basis
$(b^K_1)^{*r_1} * (b^K_2)^{*r_2} * \cdots * (b^K_m)^{*r_m}$.
Hence there is an $\eta_r \in H^{2r}(\E_2')$ such that
$(c^E_1)^*(\eta_r) \equiv c^K_r \mod{(c^K_1,\ldots,c^K_{r-1})}$.
This proves $(c_1^E)^*$ is surjective.

The Chern subring $Ch(G)$ can be generated by elements of the form $B\rho^*(x)$ 
for some representations $\rho\colon G\rightarrow U$ and homogeneous elements 
$x\in H^*(BU)$. By the first part of the proposition for any such $x$ there 
is a homogeneous $y\in H^*(\E_2')$ with $(c_1^E)^*(y)=x$. 
Hence $B\rho^*(x)=(c_1^EB\rho)^*(y)$ lies in $Ch_E(G)$.
\end{proof}

\begin{remark}{\em
For the theory $\widehat{E(n)}$, in terms of which our main theorem is
stated, we could have appealed to~\cite[3.11]{hun} instead of \cite{BH}.
}\end{remark}

\section{The category $\Ae$}
\label{section:thecat}
In this section, under the restriction to finite groups $G$, we identify
the category $\Ae$ for theories $E$ satisfying certain conditions. The
identification of $\Ae$ will depend on a positive natural number $n$
associated to $E$, and thus we begin by defining a family of categories
of elementary abelian subgroups of $G$, indexed by such $n$.

\begin{defn}
For\/ $0 \leqslant n \leqslant \infty$, define~$\An=\An(G)$ to be the
category with objects the elementary abelian $p$-subgroups of~$G$ and
morphisms the injective group homomorphisms $f \colon W \rightarrow V$
such that
$$
\forall w_1, \ldots, w_n \in W \quad \exists g \in G  \quad
\forall 1 \leqslant i \leqslant n \quad f(w_i) = g w_i g^{-1} \, .
$$
\end{defn}

\noindent
In particular, if $t$ is the $p$-rank of $G$, there are equivalences of
categories
$$\An[t]=\An[t+1]=\cdots=\An[\infty]$$
and this common category is the Quillen category~$\A$.
As noted in Example~\ref{example:A1}, it is proved in~\cite[7.1]{VCh}
that $\An[1]$~is the category~$\CR[\Ch(G)]$ of elementary abelians
associated to the Chern subring. Moreover, by~\cite[9.2]{VCh}, for each
$2 \leqslant n < \infty$ there is a subring~$R$ which satisfies the conditions
of Theorem~\ref{theorem:VCh} and has category $\CR = \An$.  The current paper
arose from the desire to find a topological construction of such a ring~$R$.

One way to explain why $\An[1]$~is the right category for~$\ECh[K](G)$
is via group characters.  On the one hand, $K^0(BG)$ is a completed
ring of (virtual) characters.  On the other hand, the morphisms in~$\An[1]$
are the group homomorphisms which preserve the values of characters for~$G$.
That is, $f \colon W \rightarrow V$ lies in~$\An[1]$ if and only if it
satisfies the equation
\begin{equation*}
\label{eqn:chars}
f^* \Res_V(\chi) = \Res_W(\chi)
\end{equation*}
for every character $\chi$~of $G$.

Switching attention to~$\An$ we recall the work of Hopkins, Kuhn and
Ravenel \cite{HKR}. The morphism $f \colon W \rightarrow V$ lies
in~$\An$ if and only if it satisfies the analagous equation for every
{\em generalised\/} character (class function) $\chi$ of $G$, {\em i.e.},
for every function on the set $G_{n,p}$ of commuting $n$-tuples of elements
of $G$ having $p$-power order which is constant on conjugacy classes and
takes values in a certain $E^*$-algebra $L(E^*)$, which is, roughly speaking,
the smallest $E^*$-algebra which contains all roots of each equation of the
form $[p^k](x) =0$: let $E^*_{cont}(B\Z_p^n) = \colim_r E^*(B(\Z/p^r)^n)$
and $S$ the multiplicatively closed subset generated by Euler classes of
continuous homomorphisms $\Z_p^n \to S^1$, then
$L(E^*) = S^{-1}E^*_{cont}(B\Z_p^n)$.

For suitable complex oriented $E$, an element $x\in E^*(BG)$ gives rise
to such a class function in essentially the following way: a commuting
$n$-tuple $(g_1,\ldots,g_n)$ in $G$ as above can be thought of as
a homomorphism $\alpha\colon\Z_p^n \to G$.
The value of a generalised character (class function) afforded by $x$ on 
$\alpha$ is then given by the composite
$E^*(BG) \xrightarrow{(B\alpha)^*} E^*_{cont}(B\Z_p^n)\lra L(E^*)$.
Theorem C of \cite{HKR} asserts that the character map $\chi_G$ associating
to each $x\in E^*(BG)$ the character it affords induces an isomorphism
$$L(E^*)\otimes_{E^*}E^*(BG) \xrightarrow{\chi_G} Cl_{n,p}(G;L(E^*))\,.$$
Thus $E^*(BG)$ is related to the
ring of generalised class functions in the same way as ordinary
characters to the $K$-theory of~$BG$.

Here a suitable theory means a complex oriented ring spectrum $E$,
whose coefficients $E^*$ are a complete, graded, local ring with maximal
ideal $\mathfrak{m}$, and residue characteristic $p>0$, such that the
mod~$\mathfrak{m}$ reduction of its formal group law has height $n$
and $p^{-1}E^*$ is non-zero.
These, together with Landweber exactness condition needed to apply
Theorem~\ref{theorem:large}, constitute some of our requirements on $E$.

We need however $E$ to satisfy one further property. Given a space $X$ and
any cohomology theory $E^*(-)$ represented by an $\Omega$-spectrum $\E_*$,
we have a `Hurewicz' construction
$$\mathcal{H}_E(X)\colon E^*(X)\longrightarrow\Hom(H^*(\E_*),H^*(X))$$
where $\Hom$ denotes, say, morphisms of graded $\F$ algebras. The map
$\mathcal{H}_E(X)$ is given by sending an element $\alpha\in E^r(X)$,
thought of as a homotopy class of maps $\alpha\colon X\rightarrow\E_r$,
to its corresponding cohomology homomorphism
$\alpha^*\colon H^*(\E_r)\rightarrow H^*(X)$.

\begin{property}\label{cond}
Say the theory $E$ satisfies {\em Property~\ref{cond}($n$)} if
\renewcommand{\labelenumi}{{(\alph{enumi})}}
\begin{enumerate}
\item $E$ is Landweber exact and satisfies the
hypotheses of~\cite{HKR}, namely: the coefficients $E^*$ are a complete,
graded, local ring with maximal ideal $\mathfrak{m}$ and residue
characteristic $p>0$, such that $p^{-1} E^*$ is nonzero and the
mod~$\mathfrak{m}$ reduction of the formal group law has height $n$;

\item the Hurewicz map $\mathcal{H}_E(X)$ is injective when $X=BV$,
the classifying space of an elementary abelian group of finite rank.
\end{enumerate}
\end{property}

\begin{remark}{\em
In the next section we prove that the Baker-W\"urgler theory
$\widehat{E(n)}$ satisfies Property~\ref{cond}($n$); it appears
that a similar proof works for Morava $E$ theory.

It may well be that the property (b) follows from (a); we know of no
examples satisfying (a) which do not also satisfy (b). On the other hand,
the Johnson-Wilson (incomplete) theories $E(n)$ are examples of Landweber
exact spectra for which property (b) holds, but not (a).

Property (b) is similar to a result of Lannes and Zarati (see, for example,
\cite[8.1]{Schwartz}) who prove such an injectivity property under the
additional assumption that the infinite loop space $\E_r$ concerned has
finite type cohomology. However, any theory $E$ satisfying part (a) will be
far from having finite type cohomology and we use very different methods in
the next section to establish the result for the $\widehat{E(n)}$.
}\end{remark}

\begin{thm}\label{theorem:maingen}
Let $G$ be a finite group and suppose $E$ satisfies
Property~{\em \ref{cond}($n$)}. Then there is a
homeomorphism of varieties
$$\Vark{Ch_E(G)} \to \colim_{\An} \Vark{H^*(BV)}\,.$$
\end{thm}

\noindent Theorem~\ref{theorem:main} will of course follow from this
and Theorem~\ref{injEhat}.

\medskip
To prove Theorem~\ref{theorem:maingen} we introduce a further category
$\CE=\CE(G)$ of elementary abelians, defined in terms of the theory $E$.

\begin{defn}\label{def:ce}
For a finite group $G$, let $\CE=\CE(G)$ be the category with objects the
elementary abelian $p$-subgroups of $G$ and morphisms the injective group
homomorphisms $f\colon W\to V$ such that $f^*\Res_V = \Res_W$ holds in
$E$-cohomology:

\begin{equation*}
\begin{diagram}[height=2em,width=2em,tight]
E^*(BW) &                & \lTo^{f^*} &                & E^*(BV) \\
        & \luTo_{\Res_W} &            & \ruTo_{\Res_V} &       \\
        &                &   E^*(BG)  &                &       \\
\end{diagram}
\end{equation*}
\end{defn}

We prove Theorem~\ref{theorem:maingen} in two stages, identifying
respectively $\An$ with $\CE$, and $\CE$ with $\Ae$. The former uses part~(a)
of \eqref{cond}; the latter needs part~(b).

\begin{prop}
\label{easybit}
Suppose the theory $E$ satisfies part {\em (a)} of
Property~{\em\ref{cond}($n$)}. Then $\CE = \An$ for all finite groups.
\end{prop}

\begin{proof}
Firstly recall from Theorem~C of \cite{HKR} that the kernel of the character
map consists entirely of $p$-torsion, as $L(E^*)$ is faithfully flat over
$p^{-1} E^*$. So the character map is injective for elementary abelian groups.

Suppose $f\colon W\to V$ is in $\An$ but not in $\CE$. Then there is an
$x\in E^*(BG)$ such that $y:=f^*\Res_V(x) - \Res_W(x) \ne 0$.
The injectivity of the character map implies that $W$ has a
rank (at most) $n$ subgroup $S$ such that $\Res_S(y) \ne 0$:
just take $S$ to be the subgroup generated by any commuting $n$-tuple on
which the generalised class function associated to~$y$ does not vanish.
Now let $T=f(S)$ and $h=f|_S$. Then $h\colon S\rightarrow T$ is an isomorphism,
and is induced by conjugation by some $g\in G$. Since conjugation by an
element of $G$ leaves $x$ fixed, we arrive at
$$\Res_S(y) =  h^* \Res_T(x) - \Res_S(x) = \Res_S(g^* x - x) = 0\,,$$
a contradiction.

Now suppose that $f\colon W\to V$ is in $\CE$ but not in $\An$. Then
there is an elementary abelian $S\subset W$ of rank at most $n$, such
that $h=f|_S\colon S\to T=f(S)$ is an isomorphism not induced by
conjugation in $G$. From the definition of $\CE$ it is clear that $h$
lies in $\CE$.
Let $g_1,\ldots,g_m$ be a (minimal) generating set for $S$, and
set $g_{m+1} = \ldots = g_n =1$ if necessary. Then
$(g_1,\ldots,g_n)$ and $(h(g_1),\ldots,h(g_n))$ are two non-conjugate
$n$-tuples in $G$, and hence are separated by a generalised class
function. Surjectivity of the character map gives us a
class~$x\in E^*(BG)$ with  $\Res_S(x) - h^* \Res_T(x)\ne 0$,
{\em i.e.\@} $\Res_V(x)$ and $f^* \Res_W(x)$ are distinct.
\end{proof}

\begin{remark}{\em
(1) In the definition of $\CE$, the condition that morphisms be mono\-morphisms
is redundant: suppose $f\colon W\to V$ has kernel $K$. Then
$\Res_K$ is trivial on $E^*(BG)$, but this cannot happen unless
$K=1$, as the character isomorphism gives a nontrivial class in the
image of restriction.

\noindent
(2) Instead of $\CE$ it might seem more appropriate to consider
a category $\C'_E$ consisting of all abelian subgroups and group
homomorphisms inducing commutative triangles as in \eqref{def:ce}.
The respective variant of Proposition~\ref{easybit} would still hold,
by essentially the same arguments. We have refrained from doing so
since our construction ultimately ends up in mod $p$ cohomology,
where the difference cannot be seen. This other approach would be
relevant were we using $p$-local or integral cohomology; compare
the final remarks in section~17 of Quillen's paper \cite{Quillen}.

\noindent
(3) By construction, every morphism in $\CE$ is also in $\Ae$.
Thus combining Proposition~\ref{easybit} with Theorem~\ref{theorem:large}
yields the chain of inclusions
$\A \subseteq \An \subseteq \Ae \subseteq \An[1]\,.$

}\end{remark}

\begin{prop}
\label{hardbit}
Suppose the theory $E$ is Landweber exact and satisfies part {\em (b)} of
Property~{\em\ref{cond}($n$)}. Then $\CE = \Ae$ for all finite groups.
\end{prop}

\begin{proof}
As just noted, it is immediate that $\CE\subseteq\Ae$ and it is the reverse
inclusion we must show. Equivalently, we need to show that the category $\CE$
does not change upon passing to the subrings of mod $p$ cohomology generated
by $E$-type characteristic classes.

So, assume $f\colon W\rightarrow V$ is in $\Ae$ but not in $\CE$.
Consider the diagram
\begin{equation*}
\begin{diagram}[height=2.5em,width=4em,tight]
E^*(BW)                &&     \lTo^{f^*}       &&            E^*(BV)   \\
              & \luTo   &                      & \ruTo &               \\
\dInto^{\mathcal{H}}   &&       E^*(BG)        && \dInto_{\mathcal{H}} \\
                       && \dTo^{\mathcal{H}}   &&                      \\
\Hom(H^*(\E_*),H^*(BW)) & \lTo & \VonH & \hLine^{f^*} &
  \Hom(H^*(\E_*),H^*(BV)) \\
              & \luTo   &                      & \ruTo &                \\
                       && \lower 5pt \hbox{$\Hom(H^*(\E_*),H^*(BG))$} &&\\
\end{diagram}
\end{equation*}
\smallskip
Injectivity of the two outside vertical maps is the assumption that $E$
satisfies part (b) of Property~\ref{cond}($n$). That $f$ lies in $\Ae$
implies the commutativity of the image of the top triangle in the bottom
one; that $f$ does not lie in $\CE$ means that the top triangle does not
commute. Commutativity of the \lq sides\rq\ of the prism follows from the
naturality of the $\mathcal{H}_E(X)$ construction with respect to maps of
the space $X$. A contradiction now follows by chasing round an element
$x\in E^*(BG)$ for which $\Res_W(x)-f^*\Res_V(x)\not=0$.
\end{proof}

\section{Injectivity of $\mathcal{H}_{\widehat{E(n)}}(BV)$}
\label{section:theproof}
The goal of this section is to prove that the $I_n$-adically complete theory
$\widehat{E(n)}$ \cite{bwa} satisfies Property~\ref{cond}($n$).
Satisfaction of part (a) is well established (we see that it is Landweber
exact in \cite{bwa} and it is noted in \cite{HKR} that the other properties
listed in part (a) also hold). Thus we must demonstrate that the maps
$$\mathcal{H}_{\widehat{E(n)}}(BV)\colon \widehat{E(n)}\null^*(BV)\rightarrow
\Hom(H^*(\EHAT_*),H^*(BV))$$
are injective for all finite rank elementary abelian $p$-groups $V$.

In fact, we shall prove the following, equivalent result in homology.

\begin{thm}\label{injEhat}
Let $V$ be a finite rank elementary abelian $p$-group. Then
$$\mathcal{H}_{\widehat{E(n)}}(BV)\colon \widehat{E(n)}\null^*(BV)\rightarrow
\Hom(H_*(BV),H_*(\EHAT_*))$$
defined by
$\big(\alpha\colon BV\rightarrow \EHAT_r\big)\mapsto\big(\alpha_* \colon
H_*(BV)\rightarrow H_*(\EHAT_r)\big)$ is injective. Here $\Hom$ denotes the
morphisms in the category of graded cocommutative $\F$ coalgebras.
\end{thm}

We start by noting that the map $\mathcal{H}_E(X)$ (in either homology or
cohomology, but from now we shall work with the homology variant) satisfies
good algebraic properties.

\begin{prop}\label{coalgalg}
Suppose that $E$ is a ring spectrum and $X$ is any space. Then
$\Hom(H_*(X),H_*(\E_*))$ has a natural $E_*$-algebra structure and
$$\mathcal{H}_E(X)\colon E^*(X)\longrightarrow \Hom(H_*(X),H_*(\E_*))$$
is a graded $E_*$-algebra homomorphism.
\end{prop}

\begin{proof}
This is essentially formal, but it is a good opportunity to introduce the
notation and operations from coalgebraic algebra that are needed in the main
proof of Theorem~\ref{injEhat} together with explicit formul\ae. For basic
references to the algebraic properties and rules of manipulation in coalgebraic
algebras (Hopf rings with further structure), see \cite{HT1,RW}.

As each $\E_r$ is an H-space, each $H_*(\E_r)$ is a Hopf algebra with
product $*$ as introduced in section~\ref{section:VCh} and coproduct
$$\psi\colon H_*(\E_r)\rightarrow H_*(\E_r) \otimes H_*(\E_r)\,.$$
However, as $E$ is a ring spectrum the graded product in $E^*(X)$ is
represented by maps $\E_r\times\E_s\rightarrow\E_{r+s}$ giving a further product
$$\circ\colon H_m(\E_r) \otimes H_n(\E_s) \rightarrow H_{m+n} (\E_{r+s})\,.$$
An element $e\in E_{-r}=E^r=\pi_0(\E_r)$, thought of as a map from a point
into $\E_r$, gives rise to an element $[e]\in H_0(\E_r)$ as the image of
$1\in H_0($point$)$. Such an element is grouplike and satisfies $[d]*[e]=[d+e]$
(where defined) and $[d]\circ [e]=[de]$; the subobject of all such elements,
written $\F[p][E^*]$, forms a sub-coalgebraic ring and $H_*(\E_*)$ is a
coalgebraic algebra over this coalgebraic ring. In fact,
$H_0(\E_*)=\F[p][E^*]$ and this should be thought of as the classical
group-ring construction, endowed with extra structure. Note that
$[0]=b_0$ and is the $*$ unit (and is distinct from 0) and $[1]$ is the
$\circ$ unit (and is distinct from 1, which, however, is identical to $[0]$).
Note also that $a\circ[0]$ is $[0]$ if $a\in H_0(\E_*)$ but is 0 otherwise.

It is of course entirely formal that the set of coalgebra maps from
$H_*(X)$, an $\F$ coalgebra, to $H_*(\E_*)$, an algebra in the category
of $\F$ coalgebras, carries itself an algebra structure. However, it will
be useful to identify the operations explicitly. Addition of say
$f,g\in\Hom(H_*(X),H_*(\E_*))$ is given by the composite
$$H_*(X)\buildrel\Delta_*\over\longrightarrow H_*(X)\otimes
H_*(X)\buildrel{f\otimes g}\over\longrightarrow H_*(\E_*)\otimes
H_*(\E_*)\buildrel *\over\longrightarrow H_*(\E_*)$$
and the product is described similarly using the $\circ$ product.
The zero element is given by the composite
$$H_*(X)\rightarrow H_*(\mbox{point})\rightarrow H_*(\E_*)$$
where the second map is that which in $H_0(-)$ sends $1$ to $[0]$;
the unit is similar, using the map representing $[1]$. Finally, the
$E_*$ action is given as follows: if $e\in E^*=E_{-*}$ and
$f\colon H_*(X)\rightarrow H_*(\E_*)$, then $ef$ is the map
$H_*(X)\rightarrow H_*(\E_*)$ sending $x\in H_*(X)$ to $[e]\circ f(x)$.
It is left to the reader to check that, with these operations,
$\mathcal{H}_E(X)$ is an $E_*$-algebra homomorphism.
\end{proof}

\begin{remark}{\em
It will be useful to note that the construction $\mathcal{H}_E(X)$ is
not only natural in the space $X$ (as used in the previous section),
but is also natural in the spectrum $E$. The strategy of the proof of
Theorem~\ref{injEhat} will be to prove the analagous result for
$\mathcal{H}_{K(n)}(B\Z/p)$, {\em i.e.} in Morava $K$-theory for
the rank 1 case, and then deduce \eqref{injEhat} from naturality in the
spectrum via the Baker-W\"urgler tower \cite{bwb} linking $K(n)$ and $\Ehat$,
and the application of an appropriate K\"unneth theorem.
}\end{remark}

Recall that the group monomorphism $\Z/p\rightarrow S^1$ induces a
homomorphism $E^*(\CP)\rightarrow E^*(B\Z/p)$. For a complex oriented
theory $E$ we shall just write $x\in E^2(B\Z/p)$ for the image of the
complex orientation $x^E$.

\begin{prop}\label{EBZp}\cite{bwb, HKR, rwb}
For $E=K(n)$ or $\Ehat$, $E^*(B\Z/p)$ is the free $E^*$ module
on basis  $\{1,x,\cdots,x^{p^n-1}\}$. Moreover, for both these
theories there is a K\"unneth isomorphism
$$E^*(B(\Z/p)^r)=E^*(B\Z/p)\tensor_{E^*}\buildrel(r)\over\cdots\tensor_{E^*}
E^*(B\Z/p).$$
\qed\end{prop}

As the homomorphism $H_*(B\Z/p)\rightarrow H_*(\CP)$ is an isomorphism in
even degrees, we shall extend the notation of section~\ref{section:VCh}
and write $\beta_r\in H_{2r}(B\Z/p)$ for the corresponding elements. Note that
the coproduct $\psi\colon H_*(B\Z/p)\rightarrow H_*(B\Z/p)\otimes H_*(B\Z/p)$
acts by $\psi(\beta_r)=\sum_{i+j=r}\beta_i\otimes\beta_j$.
Note also that $\beta_0=1$.

\begin{prop}\label{calc}
{\em (a)} Suppose $r>0$. The following formul\ae\ describe some of the
action of the map $\mathcal{H}_E(B\Z/p)$ on the class $x^r\in E^{2r}(B\Z/p)$
for a complex oriented theory $E$.
$$\begin{array}{crclc}
\mathcal{H}_E(B\Z/p)(x^r)\colon&\beta_r&\mapsto &b_1^{\circ r},&\\
&\beta_t&\mapsto &0,&\mbox{ if }\,0<t<r,\\
&\beta_0=1&\mapsto &[0]=1.&\end{array}$$

\noindent{\em (b)} {\em (}Action on $x^0$.{\em )} Suppose $e\in E_*$. Then
$$\begin{array}{crclc}
\mathcal{H}_E(B\Z/p)(e)\colon&\beta_t&\mapsto &0,&\mbox{ if }\,0<t,\qquad\\
&\beta_0=1&\mapsto &[e].\quad&\end{array}$$
\end{prop}

\begin{proof}
Given the description of the map $\mathcal{H}_E(X)$ in the proof of
Proposition~\ref{coalgalg}, these are essentially straightforward
calculations in the \lq Hopf ring calculus\rq\ of \cite{RW}.
For example, $\mathcal{H}_E(B\Z/p)(x)(\beta_t)=b_t$, by definition of $b_t$.
Then $\mathcal{H}_E(B\Z/p)(x^r)(\beta_t)$ is computed by the composite
$$\beta_t\mapsto\sum_{j_1+\cdots+j_r=t}\beta_{j_1}\otimes\cdots
\otimes\beta_{j_r}
\mapsto\sum_{j_1+\cdots+j_r=t}b_{j_1}\otimes\cdots\otimes b_{j_r}
\mapsto\sum_{j_1+\cdots+j_r=t}b_{j_1}\circ\cdots\circ b_{j_r}\,.$$
If $t=0$ there is just one term in the sum -- all $j_i=0$ and
$[0]\circ\cdots\circ[0]=[0]$. If $0<t<r$ then, in every term in the sum,
at least one $j_i=0$ and at least one $j_k>0$; thus each summand contains
the element $b_{j_i}\circ b_{j_k}=[0]\circ b_{j_k}=0$ and so the whole sum
is zero. If $t=r$ then the sum has one term in which all the $j_i=1$, giving
the $b_1^{\circ r}$ of the proposition, and all other summands are zero, as
in the previous case. Part (b) follows immediately from the definitions.
\end{proof}

\begin{prop}\label{HKn}
The morphism
$$\mathcal{H}_{K(n)}(B\Z/p)\colon K(n)^*(B\Z/p)\longrightarrow
\Hom(H_*(B\Z/p),H_*(\K_*))$$
is injective.
\end{prop}

\begin{proof}
By \eqref{EBZp} a typical element of $K(n)^*(B\Z/p)$ is of the form
$\sum_{i=0}^{p^n-1}c_ix^i$ where $c_i\in K(n)^*=\F[p][v_n,v_n^{-1}]$.
In fact, it suffices to consider homogeneous elements, given the
construction of $\mathcal{H}$. We shall suppress the unit $v_n$
(and corresponding $[v_n]\in H_*({\mathbf{K(n)}}_*)$, as in \cite{wkn})
for simplicity of notation and so shall assume $c_i\in\F$. Then generally
homogeneous elements of degree $2i$ are of the simple form $c_ix^i$,
$c_i\in\F$ and $0\leqslant i<p^n$; this is not quite true if $i\equiv0$ mod
$2(p^n-1)$ in which case the general element is $c_0+c_{p^n-1}x^{p^n-1}$.

So, for $i\not\equiv0$ mod $2(p^n-1)$, it suffices to show that, for
$c_i\not=0$,
$$\mathcal{H}_{K(n)}(B\Z/p)\left(c_ix^i\right)$$
acts non-trivially on some $\beta_t$. By Proposition~\ref{calc} it sends
$\beta_i$
to $[c_i]\circ b_1^{\circ i}$; as $b_1$ is primitive and $i>0$, this is just
$c_ib_1^{\circ i}$. From \cite{wkn} we know that all the $b_1^{\circ r}\in
H_{2r}({\mathbf{K(n)}}_{2r})$ are non-trivial if $r<p^n$.

For $i\equiv0$ mod $2(p^n-1)$ we must consider the general element
$c_0+c_{p^n-1}x^{p^n-1}$. If $c_{p^n-1}=0$ then \eqref{calc}(b) shows
$\mathcal{H}_{K(n)}(B\Z/p)(c_0)$ is non-zero on $\beta_0$. Otherwise,
assuming $c_{p^n-1}\not=0$,
$$\begin{array}{rcl}
\mathcal{H}_{K(n)}(B\Z/p)(c_0+c_{p^n-1}x^{p^n-1})\colon\beta_{p^n-1}
&\mapsto& [c_0]*([c_{p^n-1}]\circ b_1^{\circ {p^n-1}})\\
&=&[c_0]*(c_{p^n-1} b_1^{\circ {p^n-1}})\,.
\end{array}$$
As $[c_0]$ is a $*$ unit (with $*$ inverse $[-c_0]$), this last
expression is non-zero and the proof is complete.
\end{proof}

It is interesting to note that in $H_*({\mathbf{K(n)}}_*)$ although
$b_1^{\circ p^n-1}\not=0$, one more $\circ$ power of $b_1$ (or even
one more suspension) kills this element. In this sense the above proof
only `just' works.

\medskip
We now recall the tower of spectra defined in \cite{bwb} (and
implicitly in \cite{bwa}). This is a tower
$$\cdots\lra E(n)/I_n^{k+1}\lra E(n)/I_n^k\lra\cdots\lra K(n)=E(n)/I_n^1$$
where the spectrum $E(n)/I_n^k$ has homotopy $E(n)_*/I_n^k$. The
homotopy limit of this tower is the Baker-W\"urgler spectrum $\Ehat$.
Unstably, ({\em i.e.}, passing to $\Omega$ spectra), this corresponds to
a tower of fibrations of the relevant spaces; the fibre of the map
$$\mathbf{(E(n)/I_n^{k+1})}_*\lra \mathbf{(E(n)/I_n^k)}_*$$
is a product of copies of $\mathbf{K(n)}_*$ indexed by a basis of
$I_n^{k+1}/I_n^k$, {\em i.e.}, by monomials in the $v_l$, $0\leqslant l<n$
(using the convention of putting $v_0=p$) of degree $k+1$.

This tower of spectra gives rise to Baker and W\"urgler's $K(n)$ Bockstein
spectral sequence \cite{bwb}. An example of this sequence is that for the
space $B\Z/p$ in which the $E_2$-page is just
$$\Ehat\null_*\tensor_{K(n)_*}K(n)^*(B\Z/p).$$
This is entirely in even dimensions and the sequence, converging to
$(\Ehat)^*(B\Z/p)$, collapses; {\em cf.\@} Proposition~\ref{EBZp}.
We shall prove the rank 1 case of Theorem~\ref{injEhat} by examining
the Hurewicz image of this spectral sequence.

\begin{thm}\label{HEn}
The morphism
$$\mathcal{H}_{\Ehat}(B\Z/p)\colon \Ehat\null^*(B\Z/p)\longrightarrow
\Hom(H_*(B\Z/p),H_*(\EHAT_*))$$
is injective.
\end{thm}

\begin{proof}
Let $0\not=\alpha\in (\Ehat)^s(B\Z/p)$ and consider it as a map
$B\Z/p\rightarrow\EHAT_s$. Either there is some integer $k\geqslant1$
for which composition of $\alpha$ with the maps in the Baker-W\"urgler
tower gives an essential map
$\tilde\alpha\colon B\Z/p\rightarrow\mathbf{(E(n)/I_n^{k+1})}_s$
but a null map to $\mathbf{(E(n)/I_n^k)}_s$, or else $\alpha$ maps to a
non-zero element of $K(n)^s(B\Z/p)$ (at the bottom of the tower).

In the former case, the map $\tilde\alpha$ lifts to an
essential map to the fibre of the map
$\mathbf{(E(n)/I_n^{k+1})}_s\lra \mathbf{(E(n)/I_n^k)}_s$,
a product  of spaces from the $\Omega$ spectrum for $K(n)$. Thus in
both this or the second case, $\alpha$ gives rise to a non-trivial map
$$a\colon B\Z/p\longrightarrow\prod_i\K_{r_i}$$
where the product is finite, indexed by the monomials, $\w_i$ say,
in the $v_l$ of degree $k$. (Note that the $r_i$ will generally differ
from the original $s$, their value depending on the dimension of the
$\w_i$ concerned.)

Consider the commutative diagram
\begin{diagram}[height=2.5em,width=6em,tight]
H_*(B\Z/p) & \rTo & H_*(\EHAT_s)\\
\dTo^{a_*} &      &   \dTo \\
H_*(\prod_i\K_{r_i}) & \rTo & H_*(\mathbf{(E(n)/I_n^{k+1})}_s) \,.
\end{diagram}
The top map is $\mathcal{H}_{\Ehat}(B\Z/p)(\alpha)$, the map we wish
to show to be non-trivial. This will follow by finding an element
$\beta_r\in H_*(B\Z/p)$ which passes to something non-zero in
$H_*(\mathbf{(E(n)/I_n^{k+1})}_s)$ and we do this by examining the
action of $a_*$ (via \ref{HKn}) and the bottom horizontal map which
is given by the inclusion of the fibre
$\iota\colon\prod_i\K_{r_i}\rightarrow\mathbf{(E(n)/I_n^{k+1})}_s$.

First examine the map $a\colon B\Z/p\rightarrow\prod_i\K_{r_i}$. Keeping
track of the components and the monomials $\w_i$ they correspond to, we
can write this map as a tuple $(\ldots,\w_ic_{i}x^{t_i},\ldots)$ where,
as before, we assume $c_{i}\in\F$ and $0\leqslant t_i<p^n$ (and, strictly
speaking, in the dimension congruent to 0 mod $2(p^n-1)$, components may
be of the form $\w_i(c_{0,i}+c_{p^n-1,i}x^{p^n-1})$). The composite with
the inclusion, in homology,
$$H_*(B\Z/p)\buildrel a_*\over\lra H_*(\prod\K_{r_i})\lra
H_*(\mathbf{(E(n)/I_n^{k+1})}_s)$$
sends an element $\beta_r$ to
$$\sum\displaystyle\mathop{*}_i[\w_i]\circ\iota_*\mathcal{H}_{K(n)}(B\Z/p)
\left(c_{i}x^{t_i}\right)(\beta_i')$$
where the sum is over all the terms in the iterated coproduct of $\beta_r$,
and where we write $\psi(\beta_r)=\sum\cdots\otimes\beta_i'\otimes\cdots$.
(Again, reading the longer expression $c_{0,i}+c_{p^n-1,i}x^{p^n-1}$ in the
displayed formula where necessary.)

If all the powers of $x$ in this expression are zero, so that we are just
dealing with a `constant' term, then the result is easy -- taking $\beta_0$
in the top left hand of the commutative square we map to
$\sum_i[c_i\w_i]\not=0\in H_*(\mathbf{(E(n)/I_n^{k+1})}_*)$, where
the $c_i\in\F$ (not all zero) are the coefficients of $x^0$ in each
factor in the above expression for $a$.

So suppose $r$ is the smallest positive power of $x$ which appears in
the expression for $a$ above. By \eqref{HKn}, the image of $\beta_r$
in the bottom right of the commutative square is then (up to a $*$ multiple
of $*$ invertible elements $[c_{0,j}]$)
$$\sum_i[\w_i]\circ c_ib_1^{\circ r}=\left(\sum_ic_i[\w_i]\right)\circ
b_1^{\circ r}$$
where the sum is now over only some of the indexing elements $i$ (namely
those for which $t_i=r$).

It suffices now to show that expressions of the form
$\left(\sum_ic_i[\w_i]\right)\circ b_1^{\circ r}$ are non-zero in
$H_*(\mathbf{(E(n)/I_n^{k+1})}_*)$. This follows from the
following lemma, thus completing the present proof.
\end{proof}

\begin{lemma}
Let $\w_j$ be a set of monomials in the $v_l$, $0<l<n$, of degree $k$ and
of some fixed homotopy dimension and suppose $c_j\in\F$, not all zero. Then
$\left(\sum_ic_j[\w_j]\right)\circ b_1^{\circ q}$ are non-zero
in $H_*(\mathbf{(E(n)/I_n^{k+1})}_*)$ for all $0\leqslant q<p^n$.
\end{lemma}

\begin{proof}
As $\circ$ product with $b_1$ represents (double) suspension, one way to
prove this would be to observe that $\left(\sum_ic_j[\w_j]\right)$ was a
non-zero element of $H_0(\mathbf{(E(n)/I_n^{k+1})}_*)$ and that this suspends
to a non-zero element $\left(\sum_ic_j\w_j\right)$ of
$H_*(E(n)/I_n^{k+1})$. Then every intermediate
$\left(\sum_ic_j[\w_j]\right)\circ b_1^{\circ q}$
must be non-zero as well. Although the statement about
$\left(\sum_ic_j[\w_j]\right)\in H_0(\mathbf{(E(n)/I_n^{k+1})}_*)$ is true,
it is not true that this suspends to a non-zero element of the stable
gadget as, indeed, $H_*(E(n)/I_n^{k+1})=0$. However, this proof would work
if $H$ was replaced with $K(n)$ and the strategy of proof will be to deduce
the $H$ result from that for $K(n)$ by arguing with the Atiyah-Hirzebruch
spectral sequence.

First check the statements for $K(n)$. It is a basic fact on the homology
of $\Omega$ spectra that
$\left(\sum_ic_j[\w_j]\right)\not=0\in K(n)_0(\mathbf{(E(n)/I_n^{k+1})}_*)$
as this is just the group-ring \cite{RW}. To see the stable result it
suffices to check that the right unit
$$E(n)_*/I_n^{k+1}\lra K(n)_*(E(n)/I_n^{k+1})$$
is an inclusion on sums of monomials of degree $k$.
However, the cofibration of spectra
$$\bigvee \Sigma^\bullet K(n)\lra E(n)/I_n^{k+1}\lra E(n)/I_n^k$$
gives rise to a long exact sequence in Morava $K$-theory. The elements
in question map to 0 in $K(n)_*(E(n)/I_n^k)$ but lift non-trivially in
$K(n)_*(\bigvee \Sigma^\bullet K(n))$.

Recall that in $H_*(\K_*)$ the element $b_1$ suspends by
iterated $\circ$ product with itself to the elements
$b_1^{\circ r}\not=0\in H_{2r}(\K_{2r})$. These are non-zero
until we get to $b_1^{\circ p^n-1}$ (non-zero) suspending to
$b_1^{\circ p^n-1}\circ e_1=0\in H_{2p^n-1}(\K_{2p^n-1})$.
Here we write $e_1\in H_1(\K_1)$ for the single suspension
element following the notation of \cite{RW}.
However, in $K(n)_*(\K_*)$ we have $b_1^{\circ p^n-1}\circ e_1=v_ne_1$
and suspension continues indefinitely, ultimately reaching the stable
element given by the image of $1$ under the right unit in $K(n)_*(K(n))$.

Now consider the map of Atiyah-Hirzebruch spectral sequences
$$
\renewcommand{\arraystretch}{1.7}
\begin{array}{ccc}
K(n)_*\tensor_{\F}H_*(\prod\K_{r_i})&\Longrightarrow&K(n)_*(\prod\K_{r_i})\\
\Big\downarrow& & \Big\downarrow\\
K(n)_*\tensor_{\F}H_*(\mathbf{(E(n)/I_n^{k+1})}_s)&\Longrightarrow&
K(n)_*(\mathbf{(E(n)/I_n^{k+1})}_s).
\end{array}
\renewcommand{\arraystretch}{1}
$$
Take an element
$\left(\sum_jc_j[\w_j]\right)\circ b_1^{\circ q}\in H_*(\prod\K_{r_i})$
with $0\leqslant q<p^n$. We know it is a permanent (and non-trivial) cycle
in the top spectral sequence \cite{wkn} and that regarded as an element of
$K(n)_*(\prod\K_{r_i})$ it maps non-trivially to
$K(n)_*(\mathbf{(E(n)/I_n^{k+1})}_*)$.
The only way the version of this element in $H_*(-)$ could map to zero in
$H_*(\mathbf{(E(n)/I_n^{k+1})}_*)$ would be if the $K(n)_*(-)$ version
dropped Atiyah-Hirzebruch filtration in mapping to
$K(n)_*(\mathbf{(E(n)/I_n^{k+1})}_*)$. As
$\dim([v]\circ b_1^{\circ q})\leqslant2p^n-2$ this cannot happen, for
dimensional reasons.
\end{proof}

\bigskip
Theorem~\ref{injEhat} now follows from \eqref{EBZp} and \eqref{HEn}
together with the observation that $H_*(B(\Z/p)^r)$ also satisfies
a K\"unneth isomorphism. \qed

\begin{remark}{\em
The method of proof used for Theorem~\ref{injEhat} can be adapted
to establish analagous results for other theories. One of the key
elements of our proof is the use of the Baker-W\"urgler tower, and
the recent work of Baker and Lazarev \cite{BL} now allows such towers
to be built in quite general circumstances. In particular, it would
seem that height $n$ Morava $E$-theory also satisfies Property~\ref{cond}($n$).

Baker's study \cite{bak} of the homotopy type of elliptic spectra
shows that if $\wp$ is any suitable prime ideal of $Ell_*$, the
complete spectrum $(Ell)^{\widehat{\quad}}_\wp$ splits as a wedge
of suspensions of $\widehat{E(2)}$ (see \cite{bak} for precise details
of the spectra $Ell$ and ideals $\wp$ considered). Theorem~\ref{injEhat}
thus shows that all these complete spectra also satisfy Property~\ref{cond}(2)
and thus, in some sense, identifies the subring of \lq elliptic characteristic
classes\rq\ in the mod $p$ cohomology of a finite group (even though we are
still not sure what an \lq elliptic object\rq\ actually is).
}\end{remark}

\section{Examples}
\label{section:examples}

We finish by sketching some of the calculations for the chromatic subrings
$Ch_{\widehat{E(n)}}(G)$ for the simplest non-trivial example, that of $G$
the alternating group $A_4$. We take the prime $p$ to be 2.

As the $2$-rank of~$A_4$ is two, the category~$\An[2]$ is the Quillen
category~$\A$. The skeleton of $\A$ may be represented as
\[ \begin{picture}(200,40)(-100,-20)
\put(-40,0){\line(1,0){80}}
\put(-40,0){\circle*{5}}
\put( 40,0){\circle*{5}}
\put(-40,-10){\makebox(0,0){$1$}}
\put( 40,-10){\makebox(0,0){$2$}}
\put(-40, 10){\makebox(0,0){$1$}}
\put( 40, 10){\makebox(0,0){$C_3$}}
\put(  0, 10){\makebox(0,0){$1$}}
\put(-60,-10){\makebox(0,0)[r]{\emph{Rank of elem.\@ abelian}}}
\put(-60, 10){\makebox(0,0)[r]{\emph{Automorphism group}}}
\end{picture} \]
Here, the nodes represent isomorphism classes of elementary abelians,
labelled by $2$-rank and automorphism group.  The edges represent
equivalence classes of morphisms under conjugacy, the label denoting
the stabilizer in the automorphism group of the target group.
Thus here there are $3 = |C_3 : 1|$ morphisms from a rank~$1$ to a rank~$2$
elementary abelian.

The skeleton of the category~$\An[1]$ is
\[ \begin{picture}(200,40)(-100,-20)
\put(-40,0){\line(1,0){80}}
\put(-40,0){\circle*{5}}
\put( 40,0){\circle*{5}}
\put(-40,-10){\makebox(0,0){$1$}}
\put( 40,-10){\makebox(0,0){$2$}}
\put(-40, 10){\makebox(0,0){$1$}}
\put( 40, 10){\makebox(0,0){$S_3$}}
\put(  0, 10){\makebox(0,0){$C_2$}}
\end{picture} \]

\noindent Thus $\A$ is strictly contained in $\An[1]$ and so, by
Theorem~\ref{theorem:main}, the subring $Ch_{\widehat{E(1)}}(A_4)$
is strictly contained in $Ch_{\widehat{E(2)}}(A_4)$. The following
calculations demonstrate an element in the latter not in the former
as well as showing explicitly the equivalence
$\Vark{Ch_{\widehat{E(2)}}(A_4)}=\Vark{H^*(BA_4)}$.

Let~$V$ be the Sylow $2$-subgroup of~$A_4$. The Weyl group of
$V$~in $A_4$ is the cyclic group~$C_3$, permuting the non-trivial
elements transitively. As~$V$ is abelian, $H^*(BA_4)$ is the ring
of $N_{A_4}(V)$-invariants in~$H^*(BV)$.  Let $x,y$ be the basis
for~$V^*$ dual to the basis $(1 \, 2)(3 \, 4)$, $(1 \, 3)(2 \, 4)$
for~$V$.  Then we have
$$H^*(BA_4) \cong \F[2]\lbrack x,y \rbrack^{C_3} \quad |x|=|y|=1$$
where the $C_3$-action is $x \mapsto y \mapsto x+y$.
Over~$\F[4]$ we can diagonalize this action and so calculate the invariants:
$H^*(BA_4)$ is generated by $D_1$, $D_0$~and $\eta$, where
\[ D_1 = x^2 + xy + y^2 \quad
D_0 = x^2 y + x y^2 \quad \eta = x^3 + x^2 y + y^3 \, . \]
Observe that $D_1$ and $D_0$ are Dickson invariants, and that~$\eta$
is the orbit sum of~$x^2 y$.
The natural permutation representation~$\pi$ of~$A_4$ has Chern classes
$c_2(\pi) = D_1^2$ and $c_3(\pi)=D_0^2$.  These generate the Chern subring
$Ch(A_4)=Ch_{\widehat{E(1)}}(A_4)$, and are $\widehat{E(2)}$-type classes
as well by Proposition~\ref{prop:inclusion}.

Now $\widehat{E(2)}\null^*(BV) \cong \widehat{E(2)}\null^*[[w,z]] / ([2]_F(w),
[2]_F(z))$ where $[2]_F(-)$ denotes the 2-series of the formal group law for
$\widehat{E(2)}$. Set $\theta = \Tr^{BA_4}_{BV} (w^2 z)$.
By the Mackey formula, we have
\[ \Res_{BV} \theta = w^2 z + z^2 (w +_F z) + (w +_F z)^2 w\colon
BV\rightarrow\mathbf{\widehat{E(2)}}_6\, .\]
Write $\beta(s)$ for the generating function of the $\beta_i\in
H_{2i}(B\Z/2)$, that is
$\beta(s)$ is the formal power series $\sum_{i\geqslant0}\beta_is^i\in
H_*(B\Z/2)[\![s]\!]$, and similarly regard $\beta(s)\otimes\beta(t)$ as the
corresponding generating function in $H_*(BV)=H_*(B\Z/2)\tensor_{\F[2]}
H_*(B\Z/2)$.
Likewise, write $b(s)$ for the generating function of the $b_i$ in
$H_*(\mathbf{\widehat{E(2)}}_2)$. The standard Ravenel-Wilson relation
\cite[3.6]{RW} in $H_*(\E_*)$ for any complex oriented theory $E$ now
allows us to compute
\begin{multline*}
(\Res_{BV} \theta)_* (\beta(s) \otimes \beta(t)) \\
{} = \left( b(s)^{\circ2} \circ b(t) \right) *
\left( b(t)^{\circ2} \circ b(s+t) \right) *
\left( b(s+t)^{\circ2} \circ b(s) \right).
\end{multline*}
This reduces to
\[ b_0^{\circ3} + \barb(s)^{\circ2} \circ \barb(t) +
\barb(t)^{\circ2} \circ \barb(s+t) +
\barb(s+t)^{\circ2} \circ \barb(s) \quad\mbox{mod indecomposables}
\]
where $\barb(t)$ is defined as $b(t)-b_0=\sum_{r\geqslant1}b_rt^r$ (recall
that $b_r \circ b_0 = 0$ for $r>0$ and that $b_0$ is the $*$ unit).

The $(s,t)$-degree $3$ part of this expression is
\[ b_1^{\circ3} (s^2 t + t^2 (s+t) + (s+t)^2 s) =
b_1^{\circ3} (s^3 + s^2 t + t^3) \, . \]
Note that $b_1^{\circ3}\in H_6(\mathbf{\widehat{E(2)}}_6)$ represents a
non-zero indecomposable: for example, it maps to the corresponding element
in $H_6(\mathbf{K(2)}_6)$ through the Baker-W\"urgler tower, and this
represents a non-zero indecomposable here \cite{wkn}. By duality, there is
a class $\gamma \in H^6(\EHAT[2]_6)$ which kills all decomposables in~$H_6$
and sends $b_1^{\circ3}$~to $1$. Then the $\widehat{E(2)}$-type characteristic
class $\theta^*(\gamma)$ restricts to~$H^*(BV)$ as~$\eta^2$.
Hence $Ch_{\widehat{E(2)}}(A_4)$ contains $Ch_{\widehat{E(1)}}(A_4)$ as a
subring, and the inclusion of $Ch_{\widehat{E(2)}}(A_4)$ in $H^*(BA_4)$ is
an inseparable isogeny.

\end{document}